\documentclass[titlepage,twoside,12pt]{article}
\usepackage{amssymb}
\usepackage{amsfonts}
\begin{document}

{\bf \Large Solving large classes of nonlinear systems} \\

{\bf \Large of PDEs} \\

{\bf Roumen Anguelov,~Elem\'{e}r E Rosinger} \\
Department of Mathematics \\
and Applied Mathematics \\
University of Pretoria \\
Pretoria \\
0002 South Africa \\
eerosinger@hotmail.com \\

{\bf Abstract}

It is shown that large classes of nonlinear systems of PDEs, with possibly associated initial
and/or boundary value problems, can be solved by the method of order completion. The solutions
obtained can be assimilated with Hausdorff continuous functions. The usual Navier-Stokes
equations, as well as their various modifications aiming at a realistic modelling, are
included as particular cases. The same holds for the critically important constitutive
relations in various branches of Continuum Mechanics. The solution method does not involve
functional analysis, nor various Sobolev or other spaces of distributions or generalized
functions. The general and type independent {\it existence} and {\it regularity} results
regarding solutions presented here have recently been introduced in the literature. \\

\hspace{2cm} "... provided also if need be that the notion of a solution

\hspace{2.8cm} shall be suitably extended ..." \\

\hfill cited from Hilbert's 20th Problem \\

{\bf 1. Main ideas of the order completion solution method} \\

The solution method is divided in two parts. The proof of the {\it existence} of solutions
follows the method of order completion introduced and first developed in Oberguggenberger \&
Rosinger. The proof of the {\it regularity} of solutions is a consequence of recent results
obtained in Anguelov [1], regarding the structure of the Dedekind order completion of spaces
of continuous functions ${\cal C}(X)$, where $X$ is a rather arbitrary topological space. \\
For simplicity of presentation, we shall consider single nonlinear PDEs. The extension to
systems of such nonlinear PDEs and associated initial and/or boundary value problems can -
rather surprisingly - be done easily, as seen in Oberguggenberger \& Rosinger, this is being
one of the major advantages of the order completion method. Let us therefore consider
nonlinear PDEs of the general form

\smallskip
(1.1)~~~ $ F ( x, U ( x ), ~.~.~.~ , D^p_x U ( x ), ~.~.~.~ ) ~=~ f ( x ),~~~
                                               x \in \Omega ~\subseteq~ \mathbb{R}^n $

\smallskip
with $p \in \mathbb{N}^n,~~ | p | \leq m$, where the domains $\Omega$ can be any open, not
necessarily bounded subsets of $\mathbb{R}^n$, while the orders $m \in \mathbb{N}$ of the PDEs
are arbitrary given, and the unknown functions, that is, the solutions one looks for are $U :
\Omega ~\longrightarrow~ \mathbb{R}$. \\
The {\it unprecedented generality} of these nonlinear PDEs comes, above all, from the class of
functions $F$ which define the left hand terms, and which are only assumed to be {\it jointly
continuous} in all of their arguments. The right hand terms $f$ are also required to be {\it
continuous}. \\
As seen, however, both $F$ and $f$ can have certain {\it discontinuities} as well,
Oberguggenberger \& Rosinger. \\
Regardless of the above generality of the nonlinear systems of PDEs considered, one can find
for them solutions $U$ defined on the {\it whole} of the respective domains $\Omega$. These
solutions $U$ have the {\it blanket, type independent}, or {\it universal regularity}
property that they can be assimilated with {\it Hausdorff continuous functions}. \\
It follows in this way that, when solving systems of nonlinear PDEs of the generality of those
in (1.1), one can {\it dispense with} the various customary spaces of distributions,
hyperfunctions, generalized functions, Sobolev spaces, and so on. Instead one can stay within
the realms of {\it usual functions}, more precisely, of {\it interval valued} functions, see
the Appendix for a short presentation of essentials on Hausdorff continuous functions. Also,
when proving the {\it existence} and the mentioned type of {\it regularity} of such solutions
one can dispense with methods of Functional Analysis. However, functional analytic methods can
possibly be used in order to obtain further regularity or other desirable properties of such
solutions. \\
Let us now associate with each nonlinear PDE in (1.1) the corresponding nonlinear partial
differential operator defined by the left hand side, namely

\smallskip
(1.2)~~~ $ T ( x, D ) U ( x ) ~=~ F ( x, U ( x ), ~.~.~.~ , D^p_x U ( x ), ~.~.~.~ ),~~~
                                        x \in \Omega $

\smallskip
{\it Two} facts about the nonlinear PDEs in (1.1) and the corresponding nonlinear partial
differential operators $T ( x, D )$ in (1.2) are important and immediate :

\begin{itemize}

\item The operators $T ( x, D )$ can {\it naturally} be seen as acting in the {\it classical}
context, namely

\end{itemize}

(1.3)~~~ $ T ( x, D ) ~:~ {\cal C}^m ( \Omega ) \ni U ~~\longmapsto~~ T ( x, D ) U \in
                                          {\cal C}^0 ( \Omega ) $

\smallskip
while, unfortunately on the other hand :

\begin{itemize}

\item The mappings in this natural classical context (1.3) are typically {\it not} surjective
even in the case of linear $T ( x, D )$, and they are even less so in the general nonlinear
case of (1.1), (1.2).

\end{itemize}

In other words, linear or nonlinear PDEs in (1.1) typically {\it cannot} be expected to have
{\it classical} solutions $U \in {\cal C}^m ( \Omega )$, for arbitrary continuous right hand
terms $f \in {\cal C}^0 ( \Omega )$, as illustrated by a variety of well known examples, some
of them rather simple ones, see Oberguggenberger \& Rosinger [chap. 6]. \\
Furthermore, it can often happen that nonclassical solutions do have a major applicative
interest, thus they have to be sought out {\it beyond} the confines of the classical framework
in (1.3). \\
This is, therefore, how we are led to the {\it necessity} to consider {\it generalized
solutions} $U$ for PDEs like those in (1.1), that is, solutions $U \notin
{\cal C}^m ( \Omega )$, which therefore are no longer classical. This means that the natural
classical mappings (1.3) must in certain suitable ways be {\it extended} to {\it commutative
diagrams}

\begin{math}
\setlength{\unitlength}{0.2cm}
\thicklines
\begin{picture}(60,20)

\put(10,15){${\cal C}^m ( \Omega )$}
\put(27,17){$T ( x, D )$}
\put(18,15.5){\vector(1,0){26.5}}
\put(47,15){${\cal C}^0 ( \Omega )$}
\put(0,8){$(1.4)$}
\put(12,13){\vector(0,-1){8}}
\put(13.5,8){$\subseteq$}
\put(49,13){\vector(0,-1){8}}
\put(45.7,8){$\subseteq$}
\put(10,2){$~~X$}
\put(16,2.6){\vector(1,0){29.5}}
\put(47,2){$~~Y$}
\put(29,-0.5){$\widetilde{T}$}

\end{picture}
\end{math}

\smallskip
with the generalized solutions now being found as

\smallskip
(1.5) ~~~$ U \in X \setminus {\cal C}^m ( \Omega ) $

\smallskip
instead of the classical ones $U \in {\cal C}^m ( \Omega )$ which may easily fail to exist. A
further important point is that one expects to reestablish certain kind of {\it surjectivity}
type properties typically missing in (1.3), at least such as for instance

\smallskip
(1.6)~~~ $ {\cal C}^0 ( \Omega ) ~\subseteq~ \widetilde{T} ( X ) $

\smallskip
As it turns out, when constructing extensions of (1.3) given by commutative diagrams (1.4), we
shall be interested in the following somewhat larger spaces of piecewise smooth functions. For
any integer $0 \leq l \leq \infty$, we define

\smallskip
(1.7)~~~ $ {\cal C}^l_{nd} ( \Omega ) ~=~ \left \{~ u : \Omega \rightarrow \mathbb{R}~~
                    \begin{array}{|l}
                        ~\exists~ \Gamma \subset \Omega ~\mbox{closed, nowhere dense}~ : \\
                        ~~~~ u \in {\cal C}^l ( \Omega \setminus \Gamma )
                     \end{array} ~\right \} $

\smallskip
and as an immediate strengthening of (1.3), we obviously obtain

\smallskip
(1.8)~~~ $ T ( x, D )~ {\cal C}^m_{nd} ( \Omega ) ~\subseteq ~
                                         {\cal C}^0_{nd} ( \Omega ) $

\smallskip
The solution of the nonlinear PDEs in (1.1) through the order completion method will come from
the construction of specific instances of the {\it commutative diagrams} (1.4), given by, see
(2.18), (2.27)

\begin{math}
\setlength{\unitlength}{0.2cm}
\thicklines
\begin{picture}(60,30)

\put(10,25){${\cal C}^m_{nd} ( \Omega )$}
\put(27,27){$T ( x, D )$}
\put(18,25.5){\vector(1,0){26.5}}
\put(47,25){${\cal C}^0_{nd} ( \Omega )$}
\put(0,14){$(1.9)$}
\put(12,23){\vector(0,-1){6}}
\put(49,23){\vector(0,-1){6}}
\put(10,14){${\cal M}^m_T ( \Omega )$}
\put(19,14.5){\vector(1,0){25.5}}
\put(47,14){${\cal M}^0 ( \Omega )$}
\put(29,16){$T$}
\put(12,12){\vector(0,-1){6}}
\put(49,12){\vector(0,-1){6}}
\put(10,3){${\cal M}^m_T ( \Omega )^{\#}$}
\put(47,3){${\cal M}^0 ( \Omega )^{\#}$}
\put(26.5,5){$\mbox{bijective}$}
\put(19,3.5){\vector(1,0){25.5}}
\put(30,0){$T^{\#}$}

\end{picture}
\end{math}

where, as elaborated later, the operation $(~~)^{\#}$ means the {\it Dedekind order
completion}, according to MacNeille, of the respective spaces, as well as the extension to
such completions of the respective mappings, see Oberguggenberger \& Rosinger [Appendix]. It
follows that in terms of (1.4), we have

\smallskip
$ X ~=~ {\cal M}^m_T ( \Omega )^{\#},~~~ Y ~=~ {\cal M}^0 ( \Omega )^{\#},~~~
                                                      \widetilde{T} ~=~ T^{\#} $

\smallskip
thus we shall obtain for the nonlinear PDEs in (1.1) generalized solutions

\smallskip
(1.10) ~~~$ U \in {\cal M}^m_T ( \Omega )^{\#} $

\smallskip
Furthermore, instead of the {\it surjectivity} condition (1.6), we shall at least have the
following stronger one

\smallskip
(1.11) ~~~$ {\cal C}^0_{nd} ( \Omega ) ~\subseteq~ T^{\#} ( {\cal M}^m_T ( \Omega )^{\#} ) $

\smallskip
So far about the main ideas related to the {\it existence} of solutions of general nonlinear
PDEs of the form (1.1). \\
As for the {\it regularity} of such solutions, we recall that, as shown in Oberguggenberger
\& Rosinger, one has the inclusions

\smallskip
(1.12) ~~~$ {\cal M}^0 ( \Omega )^{\#} ~\subseteq~ Mes\, ( \Omega ) $

\smallskip
where $Mes\, ( \Omega )$ denotes the set of Lebesgue measurable functions on $\Omega$. In this
way, in view of (1.9), (1.10), one can assimilate the generalized solutions $U$ of the
nonlinear PDEs in (1.1) with usual measurable functions in $Mes\, ( \Omega )$. \\
Recently, however, based on results in Anguelov [1], it was shown that instead of (1.12), one
has the much stronger property

\smallskip
(1.13) ~~~$ {\cal M}^0 ( \Omega )^{\#} ~\subseteq~ \mathbb{H}\, ( \Omega ) $

\smallskip
where $\mathbb{H}\, ( \Omega )$ denotes the set of Hausdorff continuous functions on $\Omega$.
Consequently, now one can significantly improve on the earlier regularity result, as one can
assimilate the generalized solutions $U$ of the nonlinear PDEs in (1.1) with usual functions
in $\mathbb{H}\, ( \Omega )$. \\ \\

{\bf 2. The construction of diagram (1.9)} \\

Since we solve PDEs through order completion, let us see how near we can come to satisfying
the equality in (1.1), when using inequalities. For that purpose, it is useful to consider for
each $x \in \Omega$ the following set of real numbers

\smallskip
(2.1)~~~ $ \mathbb{R}_x ~=~ \{~ F ( x, \xi_0, ~.~.~.~ , \xi_p, ~.~.~.~ ) ~~|~~
       \xi_p \in \mathbb{R}, ~\mbox{for}~ p \in \mathbb{N}^n,~ | p | \leq m ~\} $

\smallskip
Clearly, for $x \in \Omega$ fixed, $\mathbb{R}_x$ is the range in $\mathbb{R}$ of
$F(x, ~.~.~.~ )$, and since $F$ is jointly continuous in all its arguments, it follows that
$\mathbb{R}_x$ is a nonvoid interval which is bounded, half bounded, or is the whole of
$\mathbb{R}$. This latter case, which can happen often with nonlinear PDEs in (1.1), will be
easier to deal with, see (2.3) next. Clearly, in the case of non-degenerate linear PDEs in
(1.1), this latter case always happens. \\
Given now $x \in \Omega$, it is obvious that a {\it necessary} condition for the existence of
a classical smooth solution $U \in {\cal C}^m$ of (1.1) in a neighbourhood of $x$ is the
condition

\smallskip
(2.2)~~~ $ f ( x ) \in \mathbb{R}_x $

\smallskip
Consequently, for the time being, we shall make the assumption that the right hand term
functions $f$ in the nonlinear PDEs in (1.1) satisfy the somewhat stronger version of
condition (2.2) given by

\smallskip
(2.3) ~~~ $ f ( x ) \in ~\mbox{int}~ \mathbb{R}_x,~~ \mbox{for}~~ x \in \Omega $

\smallskip
Clearly, whenever we have

\smallskip
(2.4)~~~ $ \mathbb{R}_x ~=~ \mathbb{R},~~ \mbox{for}~~ x \in \Omega $

\smallskip
then (2.3) is satisfied. And as mentioned, this is the case with all nontrivial linear PDEs,
as well as with most of the nonlinear PDEs of practical interest. \\
And now the basic and rather simple {\it local approximation} result on how near we can
satisfy the equality in (1.1), when using inequalities. \\

{\bf Proposition 2.1.}

\smallskip
Given $f \in {\cal C}^0 ( \Omega )$, then

\smallskip
(2.5)~~~ $ \begin{array}{l} ~\forall~~ x_0 \in \Omega,~~ \epsilon ~>~ 0 ~~: \\
               ~\exists~~ \delta ~>~ 0,~~ P ~~\mbox{polynomial in}~~ x \in \mathbb{R}^n ~~: \\
               ~\forall~~ x \in \Omega,~~ | |\, x - x_0 \,| | ~\leq~ \delta ~~: \\
                  ~~~~~ f ( x ) - \epsilon ~\leq~ T ( x, D ) P ( x ) ~\leq~ f ( x )
            \end{array} $

\smallskip
{\bf Proof} \\
Given $x_0 \in \Omega$, then for $\epsilon > 0$ small enough, condition (2.3) yields
$\xi_p \in \mathbb{R}$, with $p \in \mathbb{N}^n,~ | p | \leq m$, such that

\smallskip
(2.6) ~~ $ F ( x_0, \xi_0, ~.~.~. , \xi_p, ~.~.~.~ ) ~=~ f ( x_0 ) - \epsilon / 2 $

\smallskip
Let us take $P$ a polynomial in $x \in \mathbb{R}^n$, which satisfies the conditions

\smallskip
$ D^p_x P ( x_0 ) ~=~ \xi_p,~~ p \in \mathbb{N}^n,~ | p | \leq m $

\smallskip
In this case from (2.6) we clearly obtain the relation

\smallskip
(2.7)~~~ $ T ( x_0, D ) P ( x_0 ) ~=~ f ( x_0 ) - \epsilon / 2 $

\smallskip
and since both $T ( x, D ) P ( x )$ and $f ( x )$ are continuous in $x \in \Omega$, the local
inequality property (2.5) follows easily from (2.7).

\smallskip
The {\it global approximation} version of the inequality property in (2.5) is given in

\smallskip
{\bf Proposition 2.2.}

\smallskip
If $f \in {\cal C}^0 ( \Omega )$, then

\smallskip
(2.8)~~~ $ \begin{array}{l} ~\forall~~ \epsilon ~>~ 0 ~~: \\
                ~\exists~~ \Gamma_\epsilon ~\subset \Omega ~~\mbox{closed, nowhere dense},
                ~~ U_\epsilon \in {\cal C}^m ( \Omega \setminus \Gamma_\epsilon ) ~~: \\
              ~~~~~ f  - \epsilon ~\leq~ T ( x, D ) U_\epsilon ~\leq~ f ~~\mbox{on}~~ \Omega
                                             \setminus \Gamma_\epsilon
            \end{array} $

\smallskip
{\bf Proof} \\
Let us take a covering of $\Omega$ of the form

\smallskip
(2.9)~~~ $ \Omega ~=~ \bigcup_{\nu \in \mathbb{N}}~ K_\nu $

\smallskip
where $K_\nu$ are compact n-dimensional intervals in $\mathbb{R}^n$, namely, $K_\nu ~=~
[ a_\nu, b_\nu ] $, with $a_\nu = (a_{\nu,1}, ~.~.~. ,a_{\nu,n}),~ b_\nu = (b_{\nu,1},
~.~.~. ,b_{\nu,n}) \in \mathbb{R}^n$. We also assume, see Forster, that the covering (2.9) is
locally finite, that is

\smallskip
(2.10)~~~ $ \begin{array}{l} ~\forall~~ x \in \Omega ~~: \\
              ~\exists~~ V_x ~\subseteq~ \Omega ~~\mbox{neighbourhood of}~~ x ~~: \\
               ~~~~ \{~ \nu \in \mathbb{N}~|~ K_\nu \cap V_x \neq \phi ~\}
               ~~\mbox{is a finite set of indices}
            \end{array} $

\smallskip
and furthermore

\smallskip
(2.11)~~~$  \mbox{the interiors of}~~ K_\nu,~ \mbox{with}~~ \nu \in \mathbb{N},
                ~~\mbox{are pairwise disjoint} $

\smallskip
Let us now take $\epsilon > 0$ arbitrary but fixed. Further, we take $\nu \in \mathbb{N}$. We
shall apply Proposition 2.1 to each $x_0 \in K_\nu$. Then we obtain $\delta_{x_0} > 0$ and a
polynomial $P_{x_0}$ such that

\smallskip
$~~~~~~ f ( x ) - \epsilon ~ \leq~  T ( x, D ) P_{x_0} ( x ) ~\leq~
                     f ( x ),~~ x \in \Omega,~~ | | x - x_0 | | ~\leq~ \delta_{x_0} $

\smallskip
But $K_\nu$ is compact, therefore

\smallskip
(2.12)~~~ $ \begin{array}{l} ~ \exists~~ \delta ~>~ 0 ~~: \\
                    ~\forall~~ x_0 \in K_\nu ~~: \\
                    ~\exists~~ P_{x_0} ~~\mbox{polynomial in}~~ x \in \mathbb{R}^n ~~: \\
                    ~\forall~~ x \in \Omega,~~ | | x - x_0 | | ~\leq~ \delta_{x_0} ~~: \\
                   ~~~~~  f ( x ) - \epsilon ~ \leq~  T ( x, D ) P_{x_0} ( x ) ~\leq~ f ( x )
                 \end{array} $

\smallskip
Now we shall subdivide $K_\nu$, which was assumed to be a compact n-dimensional interval, into
n-dimensional subintervals $I_1, ~.~.~.~ , I_\mu$, so that the diameter of each of them is
less or equal $\delta$. \\
Let us denote by $J$ a generic such n-dimensional subinterval in any of the $K_\nu$, when
$\nu \in \mathbb{N}$. If $a_J \in J$ is the center of any such n-dimensional subinterval then
(2.12) gives for $x \in ~\mbox{int}~ J$

\smallskip
$~~~~~~ f ( x ) - \epsilon ~ \leq~  T ( x, D ) P_{a_J} ( x ) ~\leq~ f ( x ) $

\smallskip
Let us now take

\smallskip
(2.13)~~~ $ \Gamma_\epsilon ~=~ \Omega ~\setminus~ \bigcup_J~ \mbox{int}~ J $

\smallskip
that is, with the union ranging over all such n-dimensional subintervals $J$. In view of ???
If we define $U_\epsilon \in {\cal C}^m ( \Omega \setminus \Gamma_\epsilon )$ by

\smallskip
$ U_\epsilon ~=~ P_{a_J} ~~~\mbox{on}~~ \Omega ~\bigcap ~\mbox{int}~ J $

\smallskip
then the proof is completed.

\smallskip
{\bf Remark 2.1}

\smallskip
1) It is easy to see that the inequalities in (2.5) and (2.8) can be replaced with the
following ones, respectively

\smallskip
(2.14)~~~ $ f ( x ) ~\leq~ T ( x, D ) P ( x ) ~\leq~ f ( x ) + \epsilon $

\smallskip
(2.15)~~~ $ f ~\leq~ T ( x, D ) U_\epsilon ~\leq~ f + \epsilon $

\smallskip
as their proofs follow after the corresponding obvious minor changes in the proofs of the
above two propositions. And these four inequalities are {\it sharper} than would respectively
be the inequalities

\smallskip
$ f ( x ) - \epsilon ~\leq~ T ( x, D ) P ( x ) ~\leq~ f ( x ) + \epsilon,~~~
f - \epsilon ~\leq~ T ( x, D ) U_\epsilon ~\leq~ f + \epsilon $

\smallskip
As we shall see not much later, we do need the sharper inequalities. Indeed, the order
completion method which we shall employ is based on MacNeille's construction, therefore, it
uses {\it Dedekind cuts}. And such cuts do need the above sharper inequalities.

\smallskip
2) In Proposition 2.2., as well as in its version corresponding to the above inequality
(2.15), we have in addition the property

\smallskip
(2.16)~~~ $ mes~ ( \Gamma_\epsilon ) ~=~ 0 $

\smallskip
where $mes$ denotes the usual Lebesgue measure. Indeed, according to (2.10), (2.11) and (2.13),
$\Gamma_\epsilon$ is a countable union of rectangular grids, each generated by a finite number
of hyperplanes. \\
Here it should be noted that the presence of the {\it closed, nowhere dense} singularity sets
$\Gamma_\epsilon$ in the {\it global} inequalities (2.8) and (2.15) will prove not to be a
hindrance. And in fact, it will lead to the classes of piecewise smooth functions in (1.7)
which prove to be convenient. \\
The presence of such closed, nowhere dense {\it singularity} sets is rather deep rooted, as it
is connected with such facts as the {\it flabbiness} of related sheaves of functions, see
Oberguggenberger \& Rosinger [chapter 7], or the global version of the classical
Cauchy-Kovalevskaia theorem on analytic nonlinear PDEs, see Oberguggenberger \& Rosinger and
the literature cited there.

\smallskip
3)  As seen from the proof of Proposition 2.2., the functions $U_\epsilon$ can in fact be
chosen as piecewise polynomials in $x \in \mathbb{R}^n$.

\smallskip
Let us now note that there is an obvious ambiguity with the piecewise smooth functions in
${\cal C}^l_{nd} ( \Omega )$ in (1.7). Indeed, given any such function $u \in {\cal C}^l_{nd}
( \Omega )$, the corresponding closed, nowhere dense set $\Gamma$ cannot be defined uniquely.
Therefore, it is convenient to factor out this ambiguity, and this can be done easily as
follows. Since ${\cal C}^0_{nd} ( \Omega )$ is the largest of these spaces of functions, we
shall do for this space the mentioned factoring out, by defining on it the {\it equivalence}
relation $u \approx v$ for any two elements $u,~ v \in {\cal C}^0_{nd} ( \Omega )$, as given
by the condition

\smallskip
(2.17)~~~ $ \begin{array}{l} ~\exists~~ \Gamma \subset \Omega ~~\mbox{closed,
                                    nowhere dense}~~ : \\
                    ~~~~~~*)~~ u,~ v \in {\cal C}^0 ( \Omega \setminus \Gamma ) \\
                    ~~~~**)~~ u~=~ v ~~\mbox{on}~~ \Omega \setminus \Gamma
                \end{array} $

\smallskip
It is easy to see that $~\approx~$ defined above is indeed an equivalence relation, since the
union of a finite number of closed and nowhere dense subsets is again closed and nowhere dense.
Now we can eliminate the mentioned ambiguity by going to the {\it quotient} space

\smallskip
(2.18)~~~ $ {\cal M}^0 ( \Omega ) ~=~ {\cal C}^0_{nd} ( \Omega ) / \approx $

\smallskip
and in view of (2.8), (2.15), we define for any two elements $G,~ H \in {\cal M}^0 ( \Omega )$
the partial order $~G \leq H~$, by

\smallskip
(2.19)~~~ $ \begin{array}{l} ~\exists~~ g \in G,~~ h \in H,~~ \Gamma \subset
                                      \Omega ~~\mbox{closed, nowhere dense}~~ : \\
                ~~~~~~~*)~~ g,~ h \in {\cal C}^0 ( \Omega \setminus \Gamma ) \\
                ~~~~~**)~~ g ~\leq~ h ~~\mbox{on}~~ \Omega \setminus \Gamma
              \end{array} $

\smallskip
Let us now denote by

\smallskip
(2.20)~~~ $ ( {\cal M}^0 ( \Omega )^{\#}, \leq ) $

\smallskip
the Dedekind order completion due to MacNeille of the partially ordered space
$( {\cal M}^0 ( \Omega ), \leq )$ which was defined in (2.18), (2.19). Then this space
${\cal M}^0 ( \Omega )^{\#}$ in (2.20) is {\it order complete}, and we have the {\it order
isomorphical embedding}, see Oberguggenberger \& Rosinger [Appendix]

\smallskip
(2.21)~~~ $ {\cal M}^0 ( \Omega ) \ni G ~\longmapsto~ < G~ ] \in
                                                              {\cal M}^0 ( \Omega )^{\#} $

\smallskip
which also {\it preserves the infima and suprema}. Moreover, in view of the MacNeille
construction, we can further extend the embedding (2.21) as follows

\smallskip
\begin{math}
\setlength{\unitlength}{0.2cm}
\thicklines
\begin{picture}(60,6)

\put(10,4){${\cal M}^0 ( \Omega )$}
\put(21.5,6.5){$\mbox{o.~i.~e.}$}
\put(18.5,4.5){\vector(1,0){13}}
\put(33.5,4){${\cal M}^0 ( \Omega )^{\#}$}
\put(0,2){$(2.22)$}
\put(43,4.5){\vector(1,0){10}}
\put(46.5,6.5){$\mbox{id}$}
\put(55,4){${\cal P} ( {\cal M}^0 ( \Omega ) )$}
\put(12,0.5){$G$}
\put(17,1){\vector(1,0){14.5}}
\put(34.5,0.5){$< G~ ]$}
\put(41,1){\vector(1,0){14.5}}
\put(57,0.5){$< G~ ]$}

\end{picture}
\end{math}

In order to obtain the full situation with respect to the range of the nonlinear partial
differential operators $T ( x, D )$, we note that ${\cal C}^0 ( \Omega ) \subseteq
{\cal C}^0_{nd} ( \Omega )$, and we have the {\it order isomorphical embedding}

\smallskip
(2.23)~~~ $ {\cal C}^0 ( \Omega ) \ni g ~\longmapsto~ G \in {\cal M}^0 ( \Omega ) $

\smallskip
where $G$ is the $~\approx$ equivalence class of $g$. Furthermore, the partial order
$~\leq~$ on ${\cal M}^0 ( \Omega )$ induces on ${\cal C}^0 ( \Omega )$ through this embedding
the usual point-wise order relation of functions, namely, $g \leq h$, if and only if
$g ( x ) \leq h ( x )$, for $x \in \Omega$. \\
Let us now recall that our main interest is the construction of the commutative diagrams (1.9).
In this regard, having constructed in (2.18) - (2.20), respectively, the spaces ${\cal M}^0
( \Omega )$ and ${\cal M}^0 ( \Omega )^{\#}$ and their partial orders, the next step is to
construct the partially ordered spaces ${\cal M}^m_T ( \Omega )$ and ${\cal M}^m_T
( \Omega )^{\#}$. For that purpose we start with (1.11), namely

\smallskip
(2.24)~~~ $ T ( x, D ) : {\cal C}^m_{nd} ( \Omega ) ~\longrightarrow~
                                                  {\cal C}^0_{nd} ( \Omega ) $

\smallskip
As mentioned, we shall solve the nonlinear PDEs in (1.1) by extending through order completion
this mapping in (2.24), and we do so by constructing the commutative diagrams in (1.9). And at
this stage we are now in the position to start doing so step by step. Let us note first that
if $u \in {\cal C}^m ( \Omega \setminus \Gamma )$, where $\Gamma \subset \Omega$ is any given
closed, nowhere dense subset, then we also have

\smallskip
(2.25)~~~ $ T ( x, D ) u \in {\cal C}^0 ( \Omega \setminus \Gamma ) $

\smallskip
This means that the singularity subsets $\Gamma$ do {\it not} increase by the application of
the nonlinear partial differential operators $T ( x, D )$. However, as before with
${\cal C}^0_{nd} ( \Omega )$, the ambiguity about associating such singularity sets to
functions in ${\cal C}^m_{nd} ( \Omega )$ remains. Therefore, in view of (2.25), we shall
define the {\it equivalence} relation $u \approx_T v$, for any two functions $u,~ v \in
{\cal C}^m_{nd} ( \Omega )$, by the condition

\smallskip
(2.26)~~~ $ T ( x, D ) u ~\approx~ T ( x, D ) v $

\smallskip
which uses the equivalence relation $~\approx~$ given in (2.17), and in addition, it also
depends on the nonlinear partial differential operator $T ( x, D )$. In fact, this equivalence
relation $~\approx_T~$ on ${\cal C}^m_{nd} ( \Omega )$ is what is called the {\it pull-back}
through the mapping $T ( x, D )$ in (2.24) of the equivalence relation $~\approx~$ on
${\cal C}^0_{nd} ( \Omega )$. \\
Let us now define the quotient space

\smallskip
(2.27)~~~ $ {\cal M}^m_T ( \Omega ) ~=~ {\cal C}^m_{nd} ( \Omega ) / \approx_T $

\smallskip
in which case the mapping (2.24) generates canonically the {\it injective} mapping

\smallskip
(2.28)~~~ $ T : {\cal M}^m_T ( \Omega ) ~\longrightarrow~ {\cal M}^0 ( \Omega ) $

\smallskip
defined by $T ( U ) = G$, where $G$ is the unique $~\approx$ equivalence class in
${\cal M}^0 ( \Omega )$ of any of the $T ( x, D ) u$, where $u$ belongs to the $~\approx_T$
equivalence class $U$ in ${\cal M}^m_T ( \Omega )$. \\
At last, we can define the partial order $~\leq_T~$ on ${\cal M}^m_T ( \Omega )$ as the {\it
pull-back} through the mapping $T$ in (2.28) of the partial order $~\leq~$ in (2.19) on
${\cal M}^0 ( \Omega )$, that is, for $U,~ V \in {\cal M}^m_T ( \Omega )$, we have
$U \leq_T V$, if and only if

\smallskip
(2.29)~~~ $ T~U ~\leq~ T~ V $

\smallskip
In this way we obtain the partially ordered set $( {\cal M}^m_T ( \Omega ), \leq_T )$ giving
the desired order structure on the domain ${\cal M}^m_T ( \Omega )$ of $T$, which is the
mapping in (2.28) that corresponds now to our nonlinear partial differential operator
$T ( x, D )$ in (1.2), (1.3), or more precisely, in (2.24). \\
It is obvious in view of (2.29) that the injective mapping $T$ in (2.28) is also an {\it
order isomorphical embedding}.

\smallskip
So far, we have in this way obtained the top commutative rectangle in (1.9). \\
Applying now to $( {\cal M}^m_T ( \Omega ), \leq_T )$ the Dedekind order completion of
MacNeille, we obtain

\smallskip
(2.30)~~~ $ ( {\cal M}^m_T ( \Omega )^{\#} , \leq_T ) $

\smallskip
which is {\it order complete}, and in addition, similar with (2.21), we also have the {\it
order isomorphical embedding}

\smallskip
(2.31)~~~ $ {\cal M}^m_T ( \Omega ) \ni U ~\longmapsto~ < U~ ] \in
                                                         {\cal M}^m_T ( \Omega )^{\#} $

\smallskip
which {\it preserves the infima and suprema}. Also, similar with (2.22), we have

\bigskip
\begin{math}
\setlength{\unitlength}{0.2cm}
\thicklines
\begin{picture}(60,6)

\put(10,4){${\cal M}^m_T ( \Omega )$}
\put(21.5,6.5){$\mbox{o.~i.~e.}$}
\put(18.5,4.5){\vector(1,0){13}}
\put(33.5,4){${\cal M}^m_T ( \Omega )^{\#}$}
\put(0,2){$(2.32)$}
\put(43,4.5){\vector(1,0){10}}
\put(46.5,6.5){$\mbox{id}$}
\put(55,4){${\cal P} ( {\cal M}^m_T ( \Omega ) )$}
\put(12,0.5){$U$}
\put(17,1){\vector(1,0){14.5}}
\put(34.5,0.5){$< U~ ]$}
\put(41,1){\vector(1,0){14.5}}
\put(57,0.5){$< U~ ]$}

\end{picture}
\end{math}

\smallskip
And now all that remains is to define $T^{\#}$ in (1.9). In view of (2.32), however, this can
be done in a standard manner following from the MacNeille order completion, see
Oberguggenberger \& Rosinger [Appendix]. Consequently, one obtains the {\it order isomorphical
embedding}

\smallskip
(2.33)~~~ $ T^{\#} : {\cal M}^m_T ( \Omega ) )^{\#} ~\longrightarrow~
                                                     {\cal M}^0 ( \Omega )^{\#} $

\smallskip
which also {\it preserves the infima and the suprema}. In more detail, we have the following
commutative diagram

\bigskip
\begin{math}
\setlength{\unitlength}{0.2cm}
\thicklines
\begin{picture}(60,17)

\put(10,15){${\cal M}^m_T ( \Omega ) \ni U$}
\put(32,17){$T$}
\put(24,15.5){\vector(1,0){17.5}}
\put(45,15){$T ( U ) \in {\cal M}^0 ( \Omega )$}
\put(0,9){$(2.34)$}
\put(20,13){\vector(0,-1){8}}
\put(47,13){\vector(0,-1){8}}
\put(5,2){${\cal M}^m_T ( \Omega )^{\#} \ni~ < U~]$}
\put(24.5,2.5){\vector(1,0){14.5}}
\put(41,2){$T^{\#} ( < U~ ]~ ) ~=~$}
\put(41,-0.8){$ ~=~ < T ( U )~ ] \in {\cal M}^0 ( \Omega )$}
\put(31,-0.5){$T^{\#}$}

\end{picture}
\end{math}

\smallskip
In this way we have indeed obtained the whole of the commutative diagram (1.9), which we shall
present now in the form seen next. Here "sur" and "inj" mean mappings which are surjective
and injective, respectively, while as before, "o. i. e." means order isomorphical embedding,
and "o. i." stands for order isomorphism. The dotted arrows "$<$ - - - - - - -" mean the
"pull-back" through which the respective structures were defined

\begin{math}
\setlength{\unitlength}{0.2cm}
\thicklines
\begin{picture}(60,44)

\put(10,39){${\cal C}^m_{nd} ( \Omega )$}
\put(27,41){$T ( x, D )$}
\put(18,39.5){\vector(1,0){26.5}}
\put(47,39){${\cal C}^0_{nd} ( \Omega )$}
\put(0,21){$(2.35)$}
\put(12,37){\vector(0,-1){12}}
\put(8,31){$\mbox{sur}$}
\put(13.5,31){$\approx_T$}
\put(20,31){$<\mbox{ -~~-~~-~~-~~-~~-~~-~~-~~-~~}$}
\put(45,31){$\approx$}
\put(51,31){$\mbox{sur}$}
\put(49,37){\vector(0,-1){12}}
\put(10,22){${\cal M}^m_T ( \Omega )$}
\put(19,22.5){\vector(1,0){25.5}}
\put(47,22){${\cal M}^0 ( \Omega )$}
\put(29,24){$T$}
\put(27.5,19.5){$\mbox{o. i. e.}$}
\put(12,20){\vector(0,-1){12}}
\put(8,14){$\mbox{inj}$}
\put(13.5,14){$\leq_T$}
\put(20,14){$<\mbox{ -~~-~~-~~-~~-~~-~~-~~-~~-~~}$}
\put(45,14){$\leq$}
\put(51,14){$\mbox{inj}$}
\put(49,20){\vector(0,-1){12}}
\put(10,5){${\cal M}^m_T ( \Omega )^{\#}$}
\put(47,5){${\cal M}^0 ( \Omega )^{\#}$}
\put(19,5.5){\vector(1,0){25.5}}
\put(22,7.5){$\mbox{o. i.,~ see (3.1) below}$}
\put(30,2){$T^{\#}$}

\end{picture}
\end{math} \\

{\bf 3. General existence result} \\

One of the typical {\it main existence results} concerning the solutions of the nonlinear
PDEs in (1.1) is presented in the following theorem, see Oberguggenberger \& Rosinger [38-64]
for a proof

\smallskip
{\bf Theorem 3.1.}

\smallskip
(3.1)~~~ $ T^{\#} ~(~ {\cal M}^m_T ( \Omega )^{\#} ~) ~=~ {\cal M}^0 ( \Omega )^{\#} $

\smallskip
This means that, given the nonlinear PDEs in (1.1), for every right hand term $f \in
{\cal M}^0 ( \Omega )^{\#}$, there exists a generalized solution $U \in {\cal M}^m_T
( \Omega )^{\#}$, satisfying the relation $T^{\#} U = f$, according to the extension in
(1.9).

\smallskip
As seen in Oberguggenberger \& Rosinger [74-93], the space ${\cal M}^0 ( \Omega )^{\#}$ in
which the right hand terms $f$ of the nonlinear PDEs in (1.1) can range - and which now are
solved by Theorem 3.1. - contains a large amount of {\it discontinuous} function on $\Omega$.
Certainly, in view of (2.18), ${\cal M}^0 ( \Omega )^{\#}$ contains all the {\it piecewise
discontinuous} functions in ${\cal C}^0_{nd} ( \Omega )$.

\smallskip
What is particularly important to note is that, in view of (3.1), a variety of linear and
nonlinear PDEs can be solved, in spite of the fact that the respective PDEs are known {\it
not} to have solutions in distributions. Among them is the celebrated 1957 Hans Lewy example,
see Obrguggenberger \& Rosinger [chap. 6, 8]. \\
In this regard, it was in Oberguggenberger \& Rosinger that this Hans Lewy example of a PDE
not solvable in distribution was nevertheless solved for the first time through the method of
order completion.

\smallskip
The {\it coherence} between the solutions obtained in (3.1) and the usual classical solutions,
whenever the nonlinear PDEs in (1.1) may have the latter, follows easily from the commutative
diagram (2.35). In other words, whenever the nonlinear PDEs in (1.1) happen to have classical
solutions $U \in {\cal C}^m ( \Omega )$, then they are also generalized solution in the sense
of (3.1).

\smallskip
Finally, it is important to note that the above existence result in (3.1) can easily be
extended to {\it systems} of nonlinear PDEs of the general form in (1.1), see Oberguggenberger
\& Rosinger [chap. 8-11]. \\ \\

{\bf 4. Initial and/or boundary value problems} \\
\hspace*{0.5cm} {\bf and constitutive relations} \\

One of the significant advantages of the order completion method in solving systems of
nonlinear PDEs of the general form in (1.1) comes from the {\it ease} initial and/or boundary
value problems associated with such equations can be solved. This is in strong
contradistinction with the variety of functional analytic methods of solution where
considerable difficulties arise related to the need to restrict distributions or generalized
functions to lower dimensional manifolds. Indeed, such operations of restriction are typically
ill-defined. \\
On the other hand, when using the method of order completion, the issue of satisfying the
initial and/or boundary values can be {\it decoupled} from the issue of the existence of
solutions. Indeed, satisfying the initial and/or boundary values can be dealt with {\it first}
and {\it separately from} the issue of proving the existence of solutions. \\
Details in this regard can be found in Oberguggenberger \& Rosinger [chap. 8, 11]. And the
reason behind that rather surprising ease the order completion method exhibits when dealing
with initial and/or boundary value problems comes from a fact seen next, in section 5. \\

In Fluid Dynamics, and in general, Continuum Mechanics, a critical role is played by {\it
constitutive equations}, see Rajagopal \& Wineman, Rajagopal, or Rajagopal \& Srinivasa. \\
The usual functional analytic methods can - if at all - deal with such constitutive equations
in no less a difficult manner than they can do with initial and boundary value conditions. A
regrettable consequence of these considerable difficulties is the failure so far of functional
analytic methods to approach in any significant, let alone, systematic manner, the issue of
such critically important constitutive equations. \\
Here again, the order completion method proves its advantage by being able to deal as well
with constitutive equations. This is simply a consequence of the fact that, as mentioned next
in section 5, the order completion method can solve equations which are far more general than
the linear or nonlinear systems of PDEs, or the constitutive equations. \\ \\

{\bf 5. An abstract existence result} \\

A better understanding of the power underlying the order completion method in solving very
large classes of equations, classes {\it far beyond} any nonlinear systems of PDEs, can be
obtained from the following rather abstract existence result, see Oberguggenberger \& Rosinger
[chap. 9]. \\
Let $X$ be any set, and let $( Y, \leq )$ be any partially ordered set which has no minimum or
maximum. Further, let

\smallskip
(5.1) ~~~$ T : X ~\longrightarrow~ Y $

\smallskip
be any given mapping. The problem we consider is to find a solution $A \in X$ for the
equation

\smallskip
(5.2) ~~~$ T ( A ) ~=~ F $

\smallskip
for any given $F \in Y$. The answer is obtained as follows. Similar with the construction of
the commutative diagrams (1.9) and (2.35), one can construct commutative diagrams

\begin{math}
\setlength{\unitlength}{0.2cm}
\thicklines
\begin{picture}(60,29)

\put(14,25){$X$}
\put(30,27){$T$}
\put(18,25.5){\vector(1,0){26.5}}
\put(47,25){$Y$}
\put(0,14){$(5.3)$}
\put(15,23){\vector(0,-1){6}}
\put(47.5,23){\vector(0,-1){6}}
\put(14,14){$X_T$}
\put(19,14.5){\vector(1,0){25.5}}
\put(47,14){$Y$}
\put(30,16){$\widetilde{T}$}
\put(15,12){\vector(0,-1){6}}
\put(47.5,12){\vector(0,-1){6}}
\put(14,3){$X_T^{\#}$}
\put(47,3){$Y^{\#}$}
\put(26.5,5){$\mbox{bijective}$}
\put(19,3.5){\vector(1,0){25.5}}
\put(30,0){$T^{\#}$}

\end{picture}
\end{math}

And then the following result on the existence of solutions holds

\smallskip
{\bf Theorem 5.1.}

\smallskip
For any given $F \in Y^{\#}$, the equation

\smallskip
(5.4) ~~~$ T^{\#} ( A ) ~=~ F $

\smallskip
has a solution $A \in X_T^{\#}$, if and only if

\smallskip
(5.5) ~~~$ \begin{array}{l}
                \sup_{~Y^{\#}}~ \{~ T^{\#} ( U ) ~~|~~ U \in X_T^{\#},~~
                                             T^{\#} ( U ) \subseteq F ~\} ~=~ \\

                ~~~~~~~~=~ \inf_{~Y^{\#}}~ \{~ T^{\#} ( V ) ~~|~~ V \in X_T^{\#},~~
                                             F \subseteq T^{\#} ( V ) ~\}
            \end{array} $ \\

The significant generality of the above existence result allows, among others, the {\it
separation} mentioned in section 4, between first satisfying the initial and/or boundary value
conditions, and then, second, proving the existence of solutions in the case of general
nonlinear systems of PDEs of the form in (1.1). Indeed, by first imposing the initial and/or
boundary values, one is in fact defining the set $X$ in (5.1). And as seen above, that can be
done without any restrictions. Subsequently, condition (5.5) is both necessary and sufficient
for the existence of a generalized solution $A \in X_T^{\#}$. \\ \\

{\bf 6. The Hausdorff continuity of solutions} \\

The major {\it novelty} in this paper is about the fact that the solutions $U \in {\cal M}^m_T
( \Omega )^{\#}$ of systems of nonlinear PDEs of type (1.1), obtained according to the
procedure in Theorem 3.1., and of its generalizations can now be assimilated with {\it
Hausdorff continuous functions} in $\mathbb{H} ( \Omega )$. \\
In fact, as seen in (A.12) in the Appendix, the mentioned solutions can be assimilated with
{\it nearly finite Hausdorff continuous functions}. \\

{\bf 7. Final comments} \\

A further advantage of the order completion method is that one is {\it not} limited to
consider in (1.9) and (2.35) only the {\it pull-back} partial order $\leq_T$ generated by the
partial differential operators $T ( x, D )$ on ${\cal M}^m_T ( \Omega )$. Indeed, as seen in
Oberguggenberger \& Rosinger [chap. 13], a large variety of other partial orders on
${\cal M}^m_T ( \Omega )$ can still secure existence theorems similar to Theorem 3.1. \\
As for the use of pull-back orders, it is important to note that there exists a certain
analogy with functional analytic methods for solving PDEs. Indeed, in such methods, the
topologies considered on the domains of the partial differential operators $T ( x, D )$ are
but pull-backs through these operators of suitable topologies on their ranges. \\
Details about such facts, and in general, about certain similarities between the order
completion method and the usual functional analytic ones in solving PDEs can be found in
Oberguggenberger \& Rosinger [chap. 12].

\smallskip
The results in this paper invite a {\it comparison} with the {\it customary perception}
regarding the solution of linear or nonlinear PDEs. Typical for that perception are the
following to recent citations.

\smallskip
The 2004 edition of the Springer Universitext book "Lectures on PDEs" by V I Arnold, starts on
page 1 with the statement :

\begin{quote}

"In contrast to ordinary differential equations, there is {\it no unified theory} of partial
differential equations. Some equations have their own theories, while others have no theory at
all. The reason for this complexity is a more complicated geometry ..." (italics added)

\end{quote}

Similarly, the 1998 edition of the book "Partial Differential Equations" by L C Evans, starts
his Examples on page 3 with the statement :

\begin{quote}
"There is no general theory known concerning the solvability of all partial differential
equations. Such a theory is {\it extremely unlikely} to exist, given the rich variety of
physical, geometric, and probabilistic phenomena which can be modelled by PDE. Instead,
research focuses on various particular partial differential equations ..." (italics added)
\end{quote}

{\bf Appendix : Definition and Properties of \\
\hspace*{2.4cm} Hausdorff-Continuous Functions} \\

The Hausdorff continuous functions are not unlike the usual real valued continuous functions.
For instance, they assume real values on a dense subset of their domain of definition and are
completely determined by the values on this subset. However, these functions may also assume
{\it interval} values on a certain subset of their domain of definition. Hence the concept of
Hausdorff continuity is formulated within the realm of {\it interval valued functions}. We
shall deal in this Appendix with functions whose values can be not only usual real numbers but
also {\it extended} real numbers, that is, elements in $\overline{\mathbb{R}} = \mathbb{R}
\cup \{ -\infty, +\infty \}$. Moreover, as mentioned it proves to be convenient to allow the
values of the functions to be not only numbers in $\overline{\mathbb{R}}$, but also {\it
closed intervals} of such numbers, namely, $[ a, b ] \subseteq \overline{\mathbb{R}}$, with
$a, b \in \overline{\mathbb{R}},~ a \leq b$.

\smallskip
Towards the end of the 19-th century, Baire brought in the concepts of {\it lower} and {\it
upper semi-continuous} functions, when dealing with nonsmooth real valued functions. And in
effect, he associated with each real valued function $f$, {\it two} other real, or extended
real valued functions $I ( f )$ and $S ( f )$, with $I ( f ) \leq f \leq S ( f )$, which
proved to be particularly helpful, see (A.6), (A.7). However, following the prevailing
mentality at the time, each of these three functions was considered separately and as being a
single valued function. \\
As it turns out on the other hand, by considering {\it interval valued} functions, such as for
instance $F ( f ) = [ I ( f ), S ( f ) ]$, one can significantly improve on the understanding
and handling of non-continuous functions. \\
The study of interval valued functions can, among others, show that the particular case of
functions which have values given by one single number is appropriate for continuous functions
only. On the other hand, non-continuous functions are much better described by suitably
associated interval valued functions. \\
Indeed, in the case of functions $f$ which are {\it not} continuous, a much better description
can be obtained by considering them given by a {\it pair} of usual point valued functions,
namely $f = [~ \underline{f}, \overline{f} ~]$, thus leading to interval valued functions,
according to $f(x) = [ \underline{f}(x), \overline{f}(x) ] \subseteq \overline{\mathbb{R}}$,
for $x \in \Omega$. And then, a natural class which replaces, and also extends, the usual
point valued continuous functions is that of {\it Hausdorff-continuous} interval valued
functions, see below Definition A1. The distinctive and {\it essential} feature of these
Hausdorff-continuous functions $f = [~ \underline{f}, \overline{f} ~]$ is a condition of {\it
minimality} with respect to the {\it gap} between $\underline{f}$ and $\overline{f}$, with the
further requirement that $\underline{f}$ be lower semi-continuous, and $\overline{f}$ be upper
semi-continuous. \\
The interest in more recent times in interval valued functions comes from a number of branches of mathematics, such as
approximation theory, Sendov, and numerical analysis, Kraemer.

\smallskip
Most of the results on interval valued functions presented in this Appendix have, however,
been developed by R Anguelov. This was done in view of the usefulness of such functions in
several branches of mathematics, see Anguelov [1,2], Anguelov \& Markov, Anguelov et.al [1,2],
Anguelov \& Minani, Anguelov \& Rosinger [1-3]. Here, owing to restriction of space, we shall
only present a minimal amount of them, needed in order to support section 6 above. For the
full details in this regard, including proofs, see Appendix 2 in Anguelov \& Rosinger [3].

\smallskip
As mentioned, the class of interval valued functions of special interest here is that of {\it
Hausdorff-continuous}, or in short, {\it H-continuous} functions. As it turns out they enjoy
the {\it minimality} property (A.8) with respect to their {\it graph completion}, see (A.7),
and that allows an effective interplay between Analysis and Topology, with the latter
involving both the domain and the range of the functions dealt with. Let

\smallskip
(A.1) \quad $ \overline{\mathbb{I\, R}} ~=~ \{~ [\underline{a},~
\overline{a}] ~~|~~ \underline{a},~ \overline{a} \in
           \overline{\mathbb{R}} ~=~ \mathbb{R} \cup \{ - \infty, + \infty \},~~
                    \underline{a}~ \leq ~\overline{a}~ \} $

\smallskip
be the set of all finite or infinite closed intervals. The functions which we consider can be
defined on arbitrary topological spaces $\Omega$. For the purposes of the nonlinear PDEs
studied in this paper, however, it will be sufficient to assume that $\Omega \subseteq
\mathbb{R}^n$ are arbitrary open subsets. Let us now consider the set of interval valued
functions

\smallskip
(A.2) \quad $ \mathbb{A} ( \Omega ) ~=~ \{~ f : \Omega ~\longrightarrow~
                                                \overline{\mathbb{I\, R}} ~\} $

\smallskip
By identifying the point $a \in \overline{\mathbb{R}}$ with the degenerate interval $[a,a] \in
\overline {\mathbb{I R}}$, we consider $\overline{\mathbb{R}}$ as a subset of
$\overline{\mathbb{I R}}$. In this way $\mathbb{A} ( \Omega )$ will contain the
set of functions with extended real values, namely

\smallskip
(A.3) \quad $ {\cal A} ( \Omega ) ~=~  \{~ f : \Omega
         ~\longrightarrow~ \overline{\mathbb{R}} ~\} ~\subseteq~ \mathbb{A}( \Omega ) $

\smallskip
We define a partial order $\leq$ on $\overline{\mathbb{I R}}$ by

\smallskip
(A.4)~~~ $ [ \underline{a},~ \overline{a} ] \leq [ \underline{b},~\overline {b} ]
                    ~~~\Longleftrightarrow~~~ \underline{a}~ \leq
                   ~\underline{b},~~ \overline{a}~ \leq ~\overline{b} $

\smallskip
Now on $\mathbb{A} ( \Omega )$ we define the partial order induced by (A.4) in the usual
point-wise way, namely, for $f,~ g \in \mathbb{A} ( \Omega )$, we have

\smallskip
(A.5)~~~ $ f \leq g ~~~\Longleftrightarrow~~~ f(x) \leq g(x),~~ x \in \Omega $

\smallskip
Clearly, when restricted to ${\cal A} ( \Omega )$, the above partial order on $\mathbb{A}
( \Omega )$ reduces to the usual one among point valued functions.

\smallskip
Let $f \in \mathbb{A} ( \Omega )$. For every $x\in\Omega$, the value of $f$ is an interval,
namely, $f ( x ) ~=~ [~ \underline{f}(x),~\overline{f}(x) ~],~~~ \mbox{with}~~
\underline{f}(x),~\overline{f}(x) \in \overline{\mathbb{R}},~ \underline{f}(x) \leq
\overline{f}(x)$. Hence, every function $f \in \mathbb{A} ( \Omega )$ can be written in the
form $f ~=~ [~ \underline{f},~ \overline{f} ~],~~ \mbox{with}~~ \underline {f},~ \overline{f}
\in {\cal A} ( \Omega ),~~\underline{f} ~\leq~ f ~\leq~ \overline{f}$, and $ f \in {\cal A}
( \Omega ) ~~~\Longleftrightarrow~~~  \underline{f} ~=~ f ~=~ \overline{f} $. \\
In the particular case of functions in ${\cal A} ( \Omega )$, that is, with extended real, but point, and not nondegenerate
interval values, a number of basic results were obtained already in Baire, see also Nicolescu
for a more recent detailed presentation. The rest of the more general results concerning
functions in $\mathbb{A} ( \Omega )$, that is, with values finite or infinite closed intervals,
were developed for the first time in the above cited works of Anguelov . The few such earlier
results were obtained in Sendov, where the particular instance of $\Omega \subseteq
\mathbb{R}$ was dealt with.

\smallskip
For $x \in \Omega$, we denote by $B_\delta ( x)$ the open ball of radius $\delta$ centered at
$x$. Let us consider any {\it dense} subset $D \subseteq \Omega$, and associate with it the
pair of mappings $I ( D, \Omega, . ),~ S ( D, \Omega, . ) : \mathbb{A}( \Omega ) \rightarrow
{\cal A}(\Omega)$, called {\it lower} and {\it upper Baire operators}, respectively, where for
every function $f \in \mathbb{A} ( \Omega )$ and $x \in \Omega$, we define

\smallskip
(A.6)~~~ $ \begin{array}{l}
       I ( D, \Omega, f ) ( x ) ~=~ \sup_{\delta > 0}~ \inf~ \{~ z \in f(y) ~~|~~
                                    y \in B_\delta ( x ) \cap D ~\} \\
       S ( D, \Omega, f ) ( x ) ~=~ \inf_{\delta > 0}~ \sup~ \{~ z\in f(y) ~~|~~
                                    y \in  B_\delta ( x ) \cap D ~\}
            \end{array} $

\smallskip
In Baire, these two operators were considered and studied in the particular case of functions
$f \in {\cal A} ( \Omega )$ and when $D = \Omega$, see also Nicolescu. In view of the main
interest here in {\it interval valued} functions $f \in \mathbb{A} ( \Omega )$, it is useful
to consider as well the following third mapping, namely, $F ( D, \Omega, f ) : \mathbb{A}
(\Omega) \rightarrow \mathbb{A} ( \Omega )$, defined for $f \in \mathbb{A}(\Omega)$ by

\smallskip
(A.7)~~~ $ F ( D, \Omega, f ) ( x ) ~=~
         [~ I ( D, \Omega, f ) ( x ),~ S ( D, \Omega, f ) ( x ) ~],~~ x \in \Omega,~ $

\smallskip
and called the {\it graph completion operator}. In case $D = \Omega$, we use the simpler
notations $I ( \Omega, \Omega, f ) = I ( f ),~ S ( \Omega, \Omega, f ) = S ( f )$ and $F
( \Omega, \Omega, f ) = F ( f )$.

\smallskip
The next definition was given in Sendov in the case of $\Omega \subseteq \mathbb{R}$, however,
it can obviously be extended to any topological space $\Omega$.

\smallskip
{\bf Definition A.1}

\smallskip
A function $f \in \mathbb{A} ( \Omega )$ is called {\it Hausdorff-continuous}, or in short,
{\it H-continuous}, if and only if for every function $g \in \mathbb{A} ( \Omega )$, we have
satisfied the following {\it minimality} condition on $f$

\smallskip
(A.8)~~~ $ g ( x ) \subseteq f ( x ),~~ x \in \Omega ~~~\Longrightarrow~~~ F ( g ) ~=~ f $

\smallskip
We shall denote by $\mathbb{H} ( \Omega )$ the set of all Hausdorff-continuous interval valued
functions on $\Omega$.

\smallskip
{\bf Definition A2}

\smallskip
A function $f \in \mathbb{H}( \Omega )$ is called {\it nearly finite}, if and only if there
exists an open and dense subset $D \subseteq \Omega$, such that

\smallskip
(A.9)~~~ $ f(x) \in \overline{\mathbb{I\, R}}~~~\mbox{is a finite interval for}~~ x \in D $

\smallskip
We denote by $\mathbb{H}_{nf} ( \Omega )$ the set of nearly finite H-continuous functions $f
\in \mathbb{A} ( \Omega )$.

\smallskip
Regarding the {\it regularity} properties of solutions of general nonlinear systems of PDEs of
the form in (1.1), a crucial role is played by the following mapping

\smallskip
(A.10)~~~ $ F_0 : {\cal C}^0_{nd} ( \Omega ) \ni u ~~\longmapsto~~
                    F ( \Omega \setminus \Gamma, \Omega, u ) \in \mathbb{H}_{nf} ( \Omega ) $

\smallskip
where we recall that, according to (1.7), for every $u \in {\cal C}^0_{nd} ( \Omega )$, there
exists a closed, nowhere dense subset $\Gamma \subset \Omega$, such that $u \in {\cal C}^0
( \Omega \setminus \Gamma )$, hence in view of (A.7), $F ( \Omega \setminus \Gamma, \Omega,
u )$ is well defined. The fact that such a $\Gamma$ need not be unique, does not affect the
above definition, see Anguelov \& Rosinger [3].

The following theorem shows that the images of two functions in $C_{nd}(\Omega)$ under the
mapping $F_0$ in (A.10) are the same, if and only if these functions are equivalent with
respect to the equivalence relation (2.17).

\smallskip
{\bf Theorem A.1}

\smallskip
Let $u,v \in {\cal C}_{nd}(\Omega)$. Then $F_0(u) ~=~ F_0(v) ~~~\Longleftrightarrow~~~
u \approx v$

\smallskip
In view of (2.18), (A.10) and the above theorem now we can define a mapping

\smallskip
(A.11)~~~ $ \widetilde{F}_{0} : {\cal M}^0(\Omega) ~\longrightarrow \mathbb{H}_{nf}(\Omega) $

\smallskip
in the following way. Let $u \in U \in {\cal M}^0(\Omega)$, then $\widetilde{F}_{0}( U ) ~=~
F ( u )$. It is easy to see that the definition of $\widetilde{F}_{0}( U )$ does not depend on the
particular representative $u \in U$ of the equivalence class $U$.

\smallskip
{\bf Theorem A.2}

\smallskip
The mapping $\widetilde{F}_{0} : {\cal M}^0(\Omega) \longrightarrow \mathbb{H}_{nf}(\Omega)$
defined in (A.11) is an {\it order isomorphic embedding} with respect to the order relation
(2.19) on ${\cal M}^0(\Omega)$ and the order relation induced by (A.5) on $\mathbb{H}_{nf}
(\Omega)$. Namely, for any $U, V \in {\cal M}^0(\Omega)$, we have

\smallskip
$~~~~~~ U \leq V ~~~\Longleftrightarrow ~~~\widetilde{F}_{0}(U) \leq \widetilde{F}_{0}(V) $

\smallskip
Finally, we also have

\smallskip
{\bf Theorem A.3}

\smallskip
The set $\mathbb{H}_{nf} ( \Omega )$ is Dedekind order complete with respect to the partial
order induced on it by (A.5).

\smallskip
Let $g \in \mathbb{H}_{nf}(\Omega)$. Then there exists a subset ${\cal G} \subseteq
{\cal M}^0(\Omega)$ such that

\smallskip
$~~~~~~ g ~=~ \sup \widetilde{F}_{0}({\cal G}) ~=~ \sup \{ \widetilde{F}_{0}(G) ~|~
                                                                G \in {\cal G} \} $

\smallskip
This theorem shows that $\mathbb{H}_{nf}(\Omega)$ is the {\it smallest} Dedekind order
complete subset of $\mathbb{H}(\Omega)$ which contains the image of ${\cal M}^0(\Omega)$ under
the order isomorphical embedding $\widetilde{F}_{0}$. Hence it is {\it order isomorphic} to
the Dedekind order completion ${\cal M}^0(\Omega)^{\#}$ of ${\cal M}^0(\Omega)$. In this way
we obtain the commutative diagram, where $\widetilde{F}_{0}{}^{\#}$ denotes the order
isomorphism from ${\cal M}^0(\Omega)^\#$ to $\mathbb{H}_{nf} (\Omega)$, namely

\begin{math}
\setlength{\unitlength}{0.1cm}
\thicklines
\begin{picture}(150,60)

\put(10,10){${\cal M}^{0}(\Omega )^{\#}$}
\put(110,10){$\mathbb{H}_{nf}(\Omega )$}
\put(60,14){$\widetilde{F}_{0}^{\#}$}
\put(25,11){\vector(1,0){80}}
\put(45,6){order isomorphism}
\put(15,28){\vector(0,-1){12}}
\put(115,28){\vector(0,-1){12}}
\put(115,23){\vector(0,1){5}}
\put(0,31){${\cal M}^{0}(\Omega )= {\cal C}_{nd}^{0}(\Omega )/\approx$}
\put(65,35){$\widetilde{F}_{0}$}
\put(42,32){\vector(1,0){63}}
\put(48,27){order isomorphical embedding}
\put(110,31){$\mathbb{H}_{nf}(\Omega )$}
\put(15,49){\vector(0,-1){12}}
\put(10,52){$C_{nd}^{0}(\Omega )$}
\put(110,52){$\mathbb{H}_{nf}(\Omega )$}
\put(60,56){$F_0$}
\put(25,53){\vector(1,0){80}}
\put(45,48){graph completion}
\put(115,49){\vector(0,-1){12}}
\put(115,44){\vector(0,1){5}}

\end{picture}
\end{math}

\smallskip
Now we can bring together the above diagram with the one in (2.35) and obtain the two
successive {\it order isomorphisms}

\smallskip
\begin{math}
\setlength{\unitlength}{0.1cm}
\thicklines
\begin{picture}(150,6)

\put(0,0){$(A.12)$}
\put(20,0){${\cal M}^m_T(\Omega) )_T^{\#}$}
\put(41,1){\vector(1,0){20}}
\put(50,3){$T^{\#}$}
\put(65,0){${\cal M}^{0}(\Omega )^{\#}$}
\put(83,1){\vector(1,0){20}}
\put(90,3){$\widetilde{F}_{0}^{\#}$}
\put(107,0){$\mathbb{H}_{nf}(\Omega )$}

\end{picture}
\end{math}

\smallskip
It follows therefore that the set ${\cal M}^m_T(\Omega) )_T^{\#}$ in which the {\it solutions}
of the general nonlinear systems of PDEs of the form in (1.1) are found is mapped by the {\it
bijection} $\widetilde{F}_{0}^{\#} \circ T^{\#}$ onto the set $\mathbb{H}_{nf}(\Omega )$ of
all {\it nearly finite Hausdorff continuous functions}. Since both these mappings are order
isomorphisms, the set $\mathcal{M}^m_T(\Omega) )_T^{\#}$ is order isomorphic with the set
$\mathbb{H}_{nf}(\Omega )$. \\
Hence, the solutions of the general nonlinear systems of PDEs of the form in (1.1), which are
obtained through the order completion method, can always be assimilated with nearly finite
Hausdorff continuous functions. \\
This is the argument supporting section 6 above.

\smallskip
For the sake of further clarification, let us turn now in the remaining part of this Appendix
to some of the issues concerning the {\it discontinuities} of Hausdorff-continuous functions.
Arbitrary interval valued functions $f = [~ \underline{f},~ \overline{f} ~] \in \mathbb{A}
( \Omega )$ can exhibit a variety of types of discontinuities, and certainly not less so, than
usual point valued functions $f \in {\cal A} ( \Omega )$ do. \\
Hausdorff-continuous functions, although generalize usual point valued continuous functions,
enjoy nevertheless a number of nontrivial {\it continuity} related properties. On the other
hand, Hausdorff-continuous functions can have quite large sets of {\it discontinuities}. This
shows that they do indeed form a larger class than the usual continuous functions, even if
they still have important similar properties. \\
The fact that large enough sets of discontinuities can be present with Hausdorff-continuous
functions allows for their use - as seen in this paper - in obtaining the {\it existence of
nonclassical solutions} for large classes of systems of nonlinear PDEs. \\
Consequently, and as mentioned, Hausdorff-continuous functions - precisely since they are {\it
not} generalized functions - can be seen as {\it setting aside} to a certain extent the
variety of distributional and other traditional generalized solutions of linear and nonlinear
PDEs which have been obtained by functional analytic methods, or by the methods of the
nonlinear algebraic theory listed by the AMS Subject Classification 2000, under 46F30. \\
Indeed, solving large classes of nonlinear PDEs through Hausdorff-continuous functions offers,
among others, the following {\it double advantage} :

\begin{itemize}

\item one can bring in a significant {\it simplification} by avoiding the variety of usual
functional analytic methods with their spaces of distributions or generalized functions, and

\item one can obtain {\it universal regularity} results for solutions of large classes of
systems of nonlinear PDEs.

\end{itemize}

Theorems A.4 and A.5 below show important properties of {\it Hausdorff-continuous} functions
related to their sets of {\it discontinuities}. For every interval valued function $f =
[~ \underline{f},~ \overline{f} ~] \in \mathbb{A} ( \Omega )$, we denote by $\Gamma ( f ) ~=~
\{~ x \in \Omega ~~|~~ \underline{f} ( x ) < \overline{f} ( x )  ~\}$, which is the set of
points $x \in \Omega$ where $f$ assumes values $f ( x ) = [~ \underline{f} ( x ) ,~
\overline{f}( x ) ~]$ that are {\it non-degenerate intervals}, and not merely points. \\
It follows that at points $x \in \Gamma ( f )$, the interval valued function $f =
[~ \underline{f},~ \overline{f} ~]$ {\it cannot} be continuous in the usual sense, since it is
not a usual point valued function. \\
In the particular case of {\it Hausdorff-continuous} functions, this fact can further be
clarified. Namely, given any point $x \in \Omega$, then

\smallskip
(A.13) \quad $ \begin{array}{l}
      x \in \Gamma ( f ) ~~~\Longleftrightarrow~~~ \underline{f} ~~\mbox{and}~~ \overline{f}
                  ~~\mbox{not continuous at}~~ x ~~~\Longleftrightarrow \\
                     \Longleftrightarrow~~~
                     \underline{f} ~~\mbox{or}~~ \overline{f}
                             ~~\mbox{not continuous at}~~ x ~~~\Longleftrightarrow~~~
                                         \underline{f} (x ) ~<~ \overline{f} ( x )
                          \end{array} $

\smallskip
And now the basic result on the discontinuities of Hausdorff-continuous functions

\smallskip
{\bf Theorem A.4}

\smallskip
Given any H-continuous function $f \in \mathbb{A} ( \Omega )$. Then $\Gamma ( f )$ is of {\it
first Baire category} in $\Omega$.

\smallskip
Let us further specify the structure of the discontinuity set $\Gamma ( f )$. For
$\epsilon > 0$, let us denote $\Gamma_\epsilon ( f ) ~=~ \{~ x \in \Omega ~~|~~ \overline{f}
( x ) - \underline{f} ( x ) \geq \epsilon ~\}$. Then clearly $\Gamma ( f ) ~=~
\bigcup_{\epsilon > 0}~ \Gamma_\epsilon ( f ) ~=~\bigcup_{n \geq 1}~ \Gamma_{1 / n} ( f )$. \\
The next theorem gives a further insight into the structure of the discontinuity set
$\Gamma ( f )$  of Hausdorff-continuous functions $f \in \mathbb{A} ( \Omega )$.

\smallskip
{\bf Theorem A.5}

\smallskip
If the function $f \in \mathbb{A} ( \Omega )$ is H-continuous, then for every $\epsilon > 0$,
the set $\Gamma_\epsilon ( f )$ is {\it closed} and {\it nowhere dense} in $\Omega$.

\smallskip
An important {\it similarity} between usual continuous, and on the other hand,
Hausdorff-continuous functions is that both of them are determined {\it uniquely} if they are
known on a {\it dense} subset of their domains of definition. This property comes in spite of
the fact that, as seen in the previous section, Hausdorff-continuous functions can have
discontinuities on sets of first Baire category, and such sets can have arbitrary large
positive Lebesgue measure, see Oxtoby. Indeed, we have

\smallskip
{\bf Theorem A.6}

\smallskip
Let $f = [~ \underline{f},~ \overline{f} ~],~ g = [~ \underline{g},~ \overline{g} ~] \in
\mathbb{A} ( \Omega )$ be two H-continuous functions, and suppose given any dense subset $D
\subseteq \Omega$. Then

\smallskip
$~~~~~~ f ( x ) ~=~ g ( x ),~~ x \in D ~~~\Longrightarrow~~ f ~=~ g ~~\mbox{on}~~ \Omega $

\smallskip
The real line $\mathbb{R}$ is Dedekind order complete, but not order complete, while the
extended real line $\overline\mathbb{R}$ is both Dedekind order complete and order complete. \\
Let us recall that a partially ordered set which is order complete will also be Dedekind order complete, but as seen
above, not necessarily the other way round as well.

\smallskip
Typically, various spaces of real valued functions encountered in Analysis are neither
Dedekind order complete, nor order complete, when considered with the natural point-wise
partial order relation. \\
However, as we can see next, this situation changes when we deal with the set $\mathbb{H}
( \Omega )$ of Hausdorff-continuous functions.

\smallskip
{\bf Theorem A.7}~(Anguelov [2])

\smallskip
The set $\mathbb{H}( \Omega )$ of Hausdorff-continuous functions is order complete when
considered with the partial order in (A.5).

\smallskip
The Dedekind order completeness of the space $\mathbb{H} ( \Omega )$ of Hausdorff-continuous
functions is a {\it nontrivial} property, in view of the various connections between the usual
continuous, and on the other hand Hausdorff-continuous functions. Indeed, the space
${\cal C} ( \Omega )$ of usual real valued continuous functions, which we have seen is
strictly contained in $\mathbb{H} ( \Omega )$, is well known {\it not} to be Dedekind order
complete. On the other hand, once $\mathbb{H} ( \Omega )$ proves to be Dedekind order complete,
its order completeness follows easily from the fact that $\overline\mathbb{R}$ is order
complete. Indeed, the smallest and largest elements in $\mathbb{H} ( \Omega )$ are
respectively the functions $\Omega \ni x \longmapsto - \infty$ and $\Omega \ni x \longmapsto
+ \infty$. \\
As is well known and shown by simple examples the spaces of real valued continuous functions
${\cal C} ( \Omega )$ are {\it not} Dedekind order complete, thus, are not order complete
either. \\
Since these spaces are partially ordered in a natural way, one can apply to them the
MacNeille version of Dedekind order completion method, see Oberguggenberger \& Rosinger
[Appendix], Luxemburg \& Zaanen, or Zaanen. \\
This however being a general construction based on Dedekind type cuts, it leaves open the
question of the {\it nature} of the elements which are added to these spaces of continuous
functions by the respective Dedekind order completion process.

\smallskip
A classical, 1950 result in this regard was obtained by Dilworth in the case of {\it bounded}
and real valued continuous functions on $\Omega$ for arbitrary completely regular topological
spaces $\Omega$. Namely, the respective Dedekind order completion is given by all the normal
upper semi-continuous functions on $\Omega$. Regarding the Dedekind order completion of spaces
${\cal C} ( \Omega )$ certain results were obtained in Mack \& Johnson.

\smallskip
As seen in Anguelov [1], the Dedekind order completion of spaces ${\cal C} ( \Omega )$ of
real valued continuous functions was for the first time effectively constructed for a large
class of topological spaces $\Omega$. In this construction Hausdorff-continuous functions and
some of their subspaces play a crucial role. Consequently, both the problem of the {\it
completion} of the spaces ${\cal C} ( \Omega )$, as well as that of the {\it structure} of the
elements which are added to these spaces ${\cal C} ( \Omega )$ through the respective
completion find a convenient solution through the use of {\it interval valued} functions.

\end{document}